\documentstyle[11pt,amssymb]{amsart}

\newtheorem{lem}{Lemma}[section]
\newtheorem{thm}{Theorem}
\newtheorem{cor}{Corollary}
\newtheorem{rem}{Remark} 

\newenvironment{proof}
{{\noindent\bf Proof:~}}{\hfill$\Box$}

\def\br#1{\left(#1\right)}
\def\bra#1{(#1)}
\def\brs#1{\left\{#1\right\}}
\def\abs#1{\left|#1\right|}

\def\qcite#1{\cite{#1}}

\def\mesh{\delta}
\def\rt{\tau}
\def\tang{\nu}
\def\un{1}
\def\Tri{\Omega'} 

\def\Prob{{\Bbb P}}
\def\doint{\oint^\mesh}

\def\bl{{\cal B}}
\def\ye{{\cal Y}}

\def\curlow{\gamma^{l.c.}}
\def\curhull{\gamma^{o.b.}}
\def\curexp{\gamma^{e.p.}}

\def\lawhull{\mu^{o.b.}}
\def\lawexp{\mu^{e.p.}}

\def\lawcolper{\mu}
\def\hull{{H}}

\def\holcur{{\cal H}}
\def\holcurhull{{\cal H}^{o.b.}}
\def\holcurexp{{\cal H}^{e.p.}}
\def\holcurb{{\cal H}_M}

\def\siglaws{{\cal B}}

\begin{document}

\title{Critical percolation in the plane.\\
I. Conformal invariance and Cardy's formula.\\
II. Continuum scaling limit.}
\author{Stanislav Smirnov}
\address{
Royal Institute of Technology\\
Department of Mathematics\\
Stockholm, S10044, Sweden
\footnote{\emph{Current address:} Section de Math\'ematiques, Universit\'e de Gen\`eve,
2-4 rue du Li\`evre, Case postale 64, 1211 Gen\`eve 4, SWITZERLAND}}
\email{stas@@math.kth.se 
\footnote{\emph{Current e-mail:} Stanislav.Smirnov@@math.unige.ch}}
\date{November 15, 2001. This is a copy of an old preprint, which I will perhaps update in the future.}
\begin{abstract}
We study scaling limits and conformal invariance 
of critical site percolation on triangular lattice.
We show that some percolation-related quantities
are harmonic conformal invariants,
and calculate their values in the scaling limit.
As a particular case we obtain conformal invariance of
the crossing probabilities and Cardy's formula.
Then we prove existence, uniqueness,
and conformal invariance of the continuum scaling limit.
\end{abstract}

\maketitle

\pagestyle{myheadings}
\markboth{Stanislav Smirnov}{Critical percolation in the plane}

\section{Introduction}\label{sec:intro}

In this paper we study critical ($p=p_c=\frac12$) site percolation on triangular lattice.
First we introduce and evaluate harmonic conformal invariants,
given by the limits of certain probabilities
as mesh of the lattice tends to zero.
As a corollary we obtain crossing probabilities
(predicted by J.~Cardy in \qcite{Cardy-92})
and show their conformal invariance (conjecture attributed to M.~Aizenman 
by R.~Langlands, Ph.~Pouliot, and Y.~Saint-Aubin in \qcite{Langlands-bams}).
Then we show the existence, uniqueness, and
conformal invariance of the continuum scaling limit.
Some other similar conformal invariants
and a different approach to the scaling limit will be discussed
in a subsequent paper.
For the general background on percolation consult the
book \qcite{Grimmett-book99},
for topics related to crossing probabilities, conformal invariance,
and scaling limits see \qcite{Langlands-bams,Aizenman-statphys,Duplantier-statphys}
and other references mentioned below.

The key property of the percolation-related conformal invariants considered
is that they depend harmonically on a parameter $z$ 
(a point inside some domain $\Omega$),
which allows to determine them uniquely from their boundary behavior, and forces them
to be conformally invariant.

To motivate search for such harmonic conformal invariants,
let us first note that one can construct some assuming 
existence and conformal invariance
of the percolation scaling limit.
For example, given a simply connected domain $\Omega$ with 
two boundary points $a$ and $b$ one can define function
$h(y,z)=h_{a,b}(y,z)$
to be the expected number of different clusters spanning from
the boundary arc $ab$ to the boundary arc $ba$ and separating $y$ from $z$
(to set up this problem rigorously one actually needs
to average numbers of ``rightmost'' and ``leftmost'' cluster boundaries
separating $y$ from $z$ and taken with appropriate signs).
Then by conjectured conformal invariance we can assume $\Omega$
to be a horizontal strip of width $1$.
Symmetry arguments and translation invariance show that
$$\begin{aligned}
h(y,z)&=h(y,w)+h(w,z)~,\\
h(y,z)&=h(y+w,z+w)~,\\
h(y,z)&=0,~{\mathrm~if~}y-z\in i{\Bbb R}~,
\end{aligned}$$
it easily follows that
\begin{equation}\label{eq:strip}
h(y,z)=c\,{\mathrm Re}(y-z)~,\end{equation}
and hence function $h(y,z)$ is harmonic in both $y$ and $z$.
So the reasoning above suggests a good candidate
for a harmonic conformal invariant
(though it does not hint at the value of $c$).

One can expect that in the discrete case similar invariants should
be (almost) discrete harmonic functions,
and one should be able to obtain a ``discrete'' proof of that.
Indeed, modification of the methods described below shows that
 the scaling limit (as mesh of the lattice
tends to zero) of the discrete version of $h$
is harmonic and satisfies (\ref{eq:strip}) with $c=\sqrt{3}/4$. 
Note that when  $y$ and $z$ are taken on the boundary of $\Omega$,
the function $h$ gives the expected
number of different percolation clusters,
crossing a conformal rectangle.

There are other similar and more complicated invariants,
which we will discuss in the next version of this paper. 
For now we will concentrate on
a single harmonic conformal invariant,
which is a complexification of crossing probabilities.
Namely we will show that the probability of a point
inside a conformal triangle to be separated by a 
percolation cluster spanning two sides from the third one
is approximately a (discrete) harmonic function.
Then boundary behavior considerations allow to reconstruct
the function uniquely and force it to be conformally invariant.
There seems to be no short motivation
for this particular function to be harmonic,
but the proof is simpler than for other
analogous invariants.
Also when the point $z$ is taken on the boundary,
it gives the crossing probabilities of conformal rectangles.
Hence Cardy's formula \qcite{Cardy-92} and conformal invariance
of crossing probabilities (in the scaling limit)
follow immediately.
Moreover, one obtains enough information to prove
existence, uniqueness, and conformal invariance of
the scaling limit.

There is a classical approach (due to Shizuo Kakutani, \qcite{Kakutani-44})
to Brownian Motion along the similar lines:
exit probabilities for Brownian Motion
started at $z$ are harmonic functions in $z$ with
easily determinable boundary values.
One can prove this by showing
the same statement for the Random Walk
(discrete Laplacian of the exit probability is trivially zero),
and then passing to a limit
(Kakutani works directly with Brownian Motion).
However, our proof of discrete harmonicity is more complicated
than the Random Walk analogue:
instead of checking that Laplacian vanishes,
we find harmonic conjugates (which turn
out to be similar conformal invariants).
This requires working with the first order derivatives only,
while Laplacian involves the second order ones.
Furthemore, instead of checking Cauchy-Riemann equations
(which can be done), we check that contour integrals vanish:
it requires even less precision,
and allows for an easier passage to a limit.

Interestingly, instead of a pair of harmonic conjugate functions,
we get a ``harmonic conjugate triple.''
It seems that $2\pi/3$ rotational symmetry enters in our paper
not because of the specific lattice we consider,
but rather manifests some symmetry laws
characteristic to (continuum) percolation.

After introducing harmonic conformal invariants,
we obtain enough information to construct the continuum
scaling limit, and show its uniqueness and conformal invariance.
There are different objects which can represent the continuum scaling
limit of percolation, see the discussion in
\qcite{Aizenman-scaling,Aizenman-statphys,Aizenman-Burchard,
Aizenman-Duplantier-Aharony}, and the references therein.
In the discrete setting one defines percolation clusters
as maximal connected subgraphs of some fixed color,
and the question is what sort of object will nicely
represent the scaling limit.

If one wants to get the full information about 
the ``percolation configuration'' in the scaling limit, there seem to be two
approaches (with equivalent results).

The first one is a variation of the straightforward way, 
when one represents a percolation  configuration as
a collection of compact connected subsets of the plane,
representing different clusters.
This has to be modified, since it happens in the scaling limit
that two parts of a percolation
cluster touch without connecting (in this particular place),
since they are separated by a curve of opposite color.
However, the modification of this approach,
suggested by Michael Aizenman in \qcite{Aizenman-statphys},
works: one should represent percolation clusters
by a collection of all curves contained inside
clusters of some fixed color.
Then percolation configuration is a collection of
all curves inside all clusters of the fixed color,
such collections will almost surely
satisfy certain ``compatibility'' relations.
One can also add for convenience the collection of all curves of the
opposite color.

Another approach is to represent a percolation configuration
as a collection of ``nested'' closed curves
-- external perimeters of clusters (of both colors).
In the discrete case such curves will be simple,
in the scaling limit they cease to be simple but remain
``non-self-traversing.''
These curves are the unique curves which correspond to crossings by both colors,
and they have a canonical orientation (depending on the color of the ``outer''
side).
Then the percolation culster is a space inside some
curve $\gamma$ which is outside of all
curves of opposite orientation lying inside $\gamma$.

The philosophy of constructing the scaling limit
is to deduce first existence of subsequential limits
by compactness arguments,
and to show that quantities we know 
determine the law of the scaling limit uniquely,
so it does not depend on chosen subsequence.
We know values of crossing probabilies in the scaling limit,
and the crossing events generate any reasonable $\sigma$-algebra
for percolation  configurations.
However it is not immediate that there is a unique law
with given probabilities of crossing events
(e.g. horizontal and vertical strips generate the Borel $\sigma$-algebra
in the plane, but knowing measures of those does not
determine a Borel measure uniquely).

Thus it is easier to build the scaling limit gradually,
so we will start with some
partial aspects of the ``full percolation configuration,''
which are also of independent interest,
and sometimes are even more convinient to work with.
One can think about these either as about scaling
limit of similar objects for discrete percolation,
or as about restriction of a full percolation configuration
to a coarser $\sigma$-algebra.

Among these objects are the outer boundary of a percolation cluster,
the external perimeter of a percolation cluster,
and a lamination of a domain (telling which points on the boundary
are connected by clusters inside the domain).


\subsubsection*{Setup and notation}
We study critical site percolation on the triangular lattice with mesh $\mesh$.
For general background, consult G.~Grimmett's  and H.~Kesten's books
\qcite{Grimmett-book99,Kesten-book}.
Vertices are colored in two colors, say blue and yellow, independently
with equal probability $p=p_c=\frac12$.
By H.~Kesten's \qcite{Kesten-12}
this probability is critical  for the site percolation on the triangular lattice,
but it is interesting to note that we do not use this fact.
By a {\it blue simple path} (in a particular percolation  configuration)
we mean a sequence $\{b_j\}$ of vertices colored blue,
with adjacent $b_j$ and $b_{j+1}$.
We say that a path goes or connects to some set,
if the dual lattice hexagon around the last 
(or the first, depending on the context)
vertex intersects this set.
We will often identify this sequence of vertices with the corresponding curve
-- i.e. the broken line $b_1b_2\dots$.

For brevity we denote by $\rt$ the cube root of unity: $\rt:=\exp\br{2\pi i/3}$,
and assume that one of the lattice directions is
parallel to the real axis.

For points $a,b$ (or prime ends) on the boundary of a simply connected
domain $\Omega$ we will call {\it ``arc $ab$,''}
the part of $\partial\Omega$
between the points $a$ and $b$
(in counterclockwise direction and including the endpoints).
It need not be an arc in the strict sense, though.
Note that our considerations do not differ whether one considers
points or prime ends.
By $[ab]$ we denote the interval joining $a$ and $b$.

We denote by $\holcur$ the space of all H\"older curves 
up to parameterization, endowed with the uniform metric
(infimum of $\sup\|\gamma_1(t)-\gamma_2(t)\|_{\infty}$ over all possible
H\"older parameterizations of two curves by the interval $[0,1]$),
and by $\siglaws$ the corresponding Borel $\sigma$-algebra.

Note that for any $M$ the set $\holcurb$ of curves admitting 
H\"older parameterization
with norm at most $M$ is compact in $\holcur$
(in the uniform topology).



\subsubsection*{Acknowledgments}
I am grateful to Lennart Carleson
for getting me interested in the subject of critical percolation,
and for his encouragement and advice during this project.
For some time I attempted to follow the approach suggested by
him and Peter Jones 
(to establish Cardy's formula for equilateral triangles first,
using local geometry of triangular lattice,
global geometry of an equilateral triangle, and a study of critical cases),
so some similar elements can be found in the present paper.
I would also like to thank M.~Benedicks, I.~Binder, L.~Carleson,
P.~Jones, K.~Johansson, N.~Makarov, O.~Schramm, M.~Sodin,
and W.~Werner for valuable comments. 

\section{Harmonic conformal invariants}\label{sec:confinv}

Consider a simply connected domain $\Omega$ and three
accessible boundary points (or prime ends), 
labeled counterclockwise $a\bra{\un},a\bra{\rt},a\bra{\rt^2}$.

If the domain $\Omega$ has a smooth boundary,
there are harmonic functions 
$$h^\Omega(a\bra{\alpha},a\bra{\rt\alpha},a\bra{\rt^2\alpha},z)
=h_\alpha(z)=h_{a\bra{\alpha}}(z)~,
~~\alpha\in\brs{\un,\rt,\rt^2}~,$$ 
which are the unique solutions of the following mixed Dirichlet-Neumann problem:
\begin{equation}\label{eq:ndp}
\left\{\begin{aligned}
h_{\alpha}=1{\mathrm~at~}a\bra{\alpha}&
~,~~~h_{\alpha}=0{\mathrm~on~the~arc~}a\bra{\rt\alpha}a\bra{\rt^2\alpha}\\
\frac{\partial}{\partial(\rt\tang)}h_{\alpha}=0
&{\mathrm~on~the~arc~}a\bra{\alpha}a\bra{\rt\alpha}\\
\frac{\partial}{\partial(-\rt^2\tang)}h_{\alpha}=0
&{\mathrm~on~the~arc~}a\bra{\rt^2\alpha}a\bra{\alpha}
\end{aligned}\right.,
\end{equation}
where $\tang$ is the counterclockwise-pointing unit tangent
to $\partial\Omega$.
Moreover, functions $h_{\un},h_{\rt},h_{\rt^2}$
in some sense form a ``harmonic conjugate triple,''
see the discussion below.

As Wendelin Werner pointed out to us, function $h_\alpha$
is also given by the probability that
Brownian Motion started at $z$ and reflected on the arcs 
$a\bra{\alpha}a\bra{\rt\alpha}$, $a\bra{\rt^2\alpha}a\bra{\alpha}$
at $frac\pi3$-angle pointing towards $a\bra{\alpha}$
hits $a\bra{\alpha}$ before the arc $a\bra{\rt\alpha}a\bra{\rt^2\alpha}$.

Problem (\ref{eq:ndp}) is conformally invariant,
so we can alternatively pose it
(and this works even if $\Omega$ has non-smooth boundary) 
by conformally transferring it from $\Omega$ to some nice domain.
The solution of the problem (\ref{eq:ndp}) takes on a particularly
nice form for an equilateral triangle $\Tri$ with 
vertices $a'\bra{\un},a'\bra{\rt},a'\bra{\rt^2}$.
Then the corresponding functions become
linear functions $h'_\alpha$, which are equal to $1$ at one of the vertices
and zero at the opposite side, i.e. 
are properly rescaled distances to the sides of the triangle. 
For a general simply connected domain $\Omega$ 
there is a unique conformal map $\varphi$ of $\Omega$
to the triangle $\Tri$, sending $a\bra{\alpha}$ to $a'\bra{\alpha}$.
Therefore we remark that
\begin{rem}
$h_\alpha$ can be defined by $h_\alpha:=h'_\alpha\circ\varphi$.
\end{rem}
It is also easy to see directly
(or by mapping to an equilateral triangle)
that for a halfplane the solution of (\ref{eq:ndp})
is a hypergeometric function.

For a triangular lattice with mesh $\mesh$, we define 
an event $Q_\alpha(z),~\alpha\in\brs{\un,\rt,\rt^2},\,z\in\Omega,$ 
as an occurrence of a blue {\it simple} path going
from the arc $a\bra{\alpha} a\bra{\rt\alpha}$ to the arc 
$a\bra{\rt^2\alpha} a\bra{\alpha}$,
and separating $z$ from the arc $a\bra{\rt\alpha}a\bra{\rt^2\alpha}$.
Then we define functions 
$$H^{\Omega,\mesh}(a\bra{\alpha},a\bra{\rt\alpha},a\bra{\rt^2\alpha},z)=
H_\alpha^\mesh(z)=H_\alpha(z)
~,~~\alpha\in\brs{\un,\rt,\rt^2}~,~~z\in\Omega~,$$
to be the probabilities of $Q_\alpha(z)$'s.
Note that $H_\alpha$'s are constant in each of the lattice triangles,
so we will mainly consider their values at the vertices
of the dual hexagonal lattice (= centers of triangles).
We emphasize again, that we consider only simple paths,
neglecting therefore the points separated by the ``dangling ends.''
Otherwise the functions would have the same boundary
values as $H_\alpha$'s but would be strictly bigger
inside the domain (and therefore not harmonic).

\begin{thm}\label{thm:triangle}
As $\mesh\to0$, functions $H^\mesh_\alpha$ converge
uniformly in $\Omega$ to functions $h_\alpha$. 
\end{thm}
The statement above is a bit unrigorous, since the functions $H^\mesh_\alpha$ are defined 
only on discrete lattice depending on $\mesh$.

Since the problem (\ref{eq:ndp}) is conformally invariant,
we conclude that
\begin{cor}
The limit of $H^\mesh_\alpha$
is a conformal invariant of the points 
$a\bra{\un},a\bra{\rt},a\bra{\rt^2},z$
and the domain $\Omega$.
\end{cor}
Consider a conformal rectangle, i.e. a simply connected domain $\Omega$
with four points (or prime ends) $a\bra{\un},a\bra{\rt},a\bra{\rt^2},x$
on the boundary, labeled counterclockwise.
Then the {\it crossing probability},
i.e. the probability of having a blue cluster connecting
the arc $xa\bra{\un}$ to the arc $a\bra{\rt}a\bra{\rt^2}$ 
is equal to $H^\mesh_{\rt^2}(x)$ 
(or rather to the boundary value of $H^\mesh_{\rt^2}$ at $x$). 
Therefore
\begin{cor}
As $\mesh\to0$, the crossing probability 
tends to $h_{\rt^2}(x)$,
and hence is conformally invariant.
\end{cor}
For a particular ``nice'' domain, like a rectangle or a half plane,
it is easy to find $h_{\alpha}$'s and then check that
crossing probabilities satisfy the Cardy's formula \qcite{Cardy-92}.
This is particularly easy for an equilateral triangle, where $h_\beta$'s are linear:
\begin{cor}[Cardy's formula in Carleson's form] 
Consider an equilateral triangle $abc$
with side length one
and a point $x\in[bc]$.
As $\mesh\to0$, the probability of a blue crossing
from $[xb]$ to $[ac]$,
tends to $|xb|$.
\end{cor}

Sometimes it is more convenient to work with a different
formulation of the results above.
For a discrete percolation configuration in conformal triangle $abc$
there is always a unique vertex $w$ on the arc $bc$
with neighbors accessible by a blue crossing from $ac$
and a yellow crossing from $ab$
(it is just the endpoint 
of a blue crossing from $ac$ to $bc$, closest to $ab$, 
or equivalently
of a yellow crossing from $ab$ to $bc$, closest to $ac$).
It follows immediately from the discussion above that
\begin{cor}[Cardy-Carleson law]
The law of $w$ converges as $\mesh\to0$
to a distribution, which is a conformal invariant
of configuration $(\Omega,a,b,c)$.
For an equilateral triangle $abc$ this distribution
is uniform on the side $bc$.
\end{cor}
\begin{rem}\label{rem:hull}
Note also, that the probability of a point $z\in\Omega$
to be contained between the two crossings tends to
$1-h_{b}(z)-h_{c}(z)\equiv h_{a}(z)$.
\end{rem}

One can also ask about the speed of convergence in the statements above.
\begin{rem}
One can check that the speed of convergence 
in the Theorem~\ref{thm:triangle} and
its Corollaries above is 
$O(\mesh^\varepsilon)$ for some $\varepsilon>0$,
and then obtain values of crossing probabilities
for the discrete case up to such error.
After establishing values of exponents
for percolation one can get $\varepsilon=2/3$.
\end{rem}

%

The remaining part of this Section is devoted to the proof 
of Theorem~\ref{thm:triangle}.

\subsection{Proof of Theorem~\ref{thm:triangle}}

Take $\beta\in\brs{1,\rt,\rt^2}$.
If $z$ is a center of some lattice triangle, and $z+\eta$
is the center of one of the adjacent triangles,
denote by $P_\beta(z,\eta)$ probability of the event
$Q_\beta(z+\eta)\setminus Q_\beta(z)$.
Then the discrete derivative of $H_\beta$
can be written as the following difference of the probabilities:
\begin{equation}\label{eq:der}
\frac{\partial}{\partial\eta}H_\beta(z)~:=~
H_\beta(z+\eta)-H_\beta(z)~=~
~P_\beta(z,\eta)-P_\beta(z+\eta,-\eta)~.
\end{equation}

\begin{rem}
Since the discrete derivative is expected to
have order $\asymp\mesh$, and
is a difference of the two probabilities on the right hand side,
one might be tempted to assume that the latter also
have order $\asymp\mesh$.
Surprisingly, they have a much bigger order $\asymp\mesh^{2/3}$,
so there is a non-trivial cancelation. 
\end{rem}

\begin{lem}[$2\pi/3$-Cauchy-Riemann equations]\label{lem:rot}
Let $z$ be a center of some triangle,
and $\eta$ be a vector from $z$ to the center of one of 
the adjacent triangles. 
Then for any $\beta\in\brs{1,\rt,\rt^2}$
\begin{equation}\label{eq:rot}
P_\beta(z,\eta)~=~P_{\rt\beta}(z,\rt\eta)~.
\end{equation}
\end{lem}

\begin{rem}
Lemma~\ref{lem:rot} holds not only for triangular lattice,
but for any graph which is a triangulation.
\end{rem}

\begin{rem}
Evaluating one more discrete derivative 
(by shifting the boundary of the domain, rather than the point),
one can check that for nearby
points $z'$ and $z''$
\begin{equation}\label{eq:c1p}
P_\beta(z',\eta')~=~P_\beta(z'',\eta')
~+~O(\mesh^{1+\varepsilon})~.
\end{equation}
Together with the Lemma~\ref{lem:rot} this gives
$$\begin{aligned}
\frac{\partial}{\partial\eta}H_\beta(z)&=
P_\beta(z,\eta)-P_\beta(z+\eta,-\eta)\stackrel{(\ref{eq:rot})}{=} 
P_{\rt\beta}(z,\rt\eta)-P_{\rt\beta}(z+\eta,-\rt\eta)\\
&\stackrel{(\ref{eq:c1p})}{=} 
P_{\rt\beta}(z,\rt\eta)-P_{\rt\beta}(z+\rt\eta,-\rt\eta)+O(\mesh^{1+\varepsilon})=
\frac{\partial}{\partial(\rt\eta)}H_{\rt\beta}(z)+O(\mesh^{1+\varepsilon})~.
\end{aligned}$$
Therefore the functions $H_{\un},H_{\rt},H_{\rt^2}$
satisfy discrete Cauchy-Riemann equations
(or rather their version for triples) up to $\mesh^\varepsilon$
(actually even up to $\mesh^{2/3}$).
But in passing to a scaling limit it is better to use
global manifestations of analyticity, so we work with
vanishing contour integrals instead.
\end{rem}

\begin{proof}
Name the vertices of the lattice triangle which contains $z$
by the letters $X,Y,Z$ starting with the one opposite to $z+\eta$ and going 
counterclockwise.
For the event $Q':=Q_\beta(z+\eta)\setminus Q_\beta(z)$ to occur,
the closest to the arc $a\bra{\rt\beta}a\bra{\rt^2\beta}$ blue
simple path $\gamma$ going from the arc $a\bra{\beta} a\bra{\rt\beta}$ 
to the arc $a\bra{\rt^2\beta} a\bra{\beta}$
should separate $z$ from $z+\eta$.
Firstly this means that
there are two disjoint blue paths (``halves'' of $\gamma$),
which go from $Y$ and $Z$ to the arcs
$a\bra{\rt^2\beta}a\bra{\beta}$ and $a\bra{\beta} a\bra{\rt\beta}$
respectively.
Secondly, the vertex $X$ is colored yellow
(otherwise we can include it into the path $\gamma$),
and is joined by a simple yellow path to
$a\bra{\rt\beta} a\bra{\rt^2\beta}$.

So we conclude that the event $Q'$
can be described as an occurrence
of three disjoint simple paths,
joining $X$, $Y$, and $Z$ to the arcs
$a\bra{\rt\beta} a\bra{\rt^2\beta}$,
$a\bra{\rt^2\beta} a\bra{\beta}$, and $a\bra{\beta} a\bra{\rt\beta}$,
and colored yellow, blue, and blue
correspondingly.

But in the latter description we can easily change the colors of the paths
while preserving the probability of $Q'$. 
In fact, for a given configuration we can chose
the ``counterclockwise-most'' yellow path from
$X$ to the arc $a\bra{\rt\beta} a\bra{\rt^2\beta}$,
and the ``clockwise-most'' blue path from
$Y$ to the arc $a\bra{\rt^2\beta} a\bra{\beta}$.
Denote by $\Omega'$ the union of these two paths
and the part of $\Omega$ between them,
containing $a\bra{\rt^2\beta}$.
The existence of a blue 
path from $Z$ to $a\bra{\beta} a\bra{\rt\beta}$ 
depends only on the coloring of $\Omega\setminus\Omega'$.
Thus if we condition by the coloring of $\Omega'$,
the probability of the existence of a blue 
path from $Z$ to $a\bra{\beta} a\bra{\rt\beta}$ 
is the same as the probability of such a yellow path,
since we can invert colors in $\Omega\setminus\Omega'$
and use that $p=\frac12$.

Taking expectation over all possible configurations of $\Omega'$,
we deduce, that the event $Q'$ has the same probability as the 
occurrence of three disjoint simple paths,
joining $X$, $Y$, and $Z$ to the arcs
$a\bra{\rt\beta} a\bra{\rt^2\beta}$,
$a\bra{\beta} a\bra{\rt^2\beta}$, and $a\bra{\beta} a\bra{\rt\beta}$,
and colored yellow, blue, and yellow
correspondingly.
But if one inverts colors in all of $\Omega$
(which preserves probabilities),
this is the description 
of the event
$Q_{\rt\beta}(z+\rt\eta)\setminus Q_{\rt\beta}(z)$,
and we proved the Lemma.
\end{proof}

\begin{lem}[H\"older norm estimates]\label{lem:hold}
There are constants $\varepsilon$ and $C$ depending
on the domain $\Omega$ only, such that
$H_\beta$ has $\varepsilon$-H\"older norm at most $C$.
The boundary values of
$H_\beta$ are zero on the arc $a\bra{\rt\beta}a\bra{\rt^2\beta}$
and tend to $1$ at the point $a\bra{\beta}$ as $\mesh\to0$.
\end{lem}

\begin{rem}
It is also easy to show that on $\partial\Omega$
$$H_{\un}+H_{\rt}+H_{\rt^2}\to1~.$$
\end{rem}

\begin{proof}
Clearly, $H_\beta$ takes values in $[0,1]$.
To prove that it is H\"older, we need to show that
\begin{equation}\label{eq:hold}
\abs{H_\beta(z)-H_\beta(z')}~\le~C'\abs{z-z'}^{\varepsilon}~,
\end{equation}
and it is sufficient to do so for pairs $z,z'$
well away from one of the boundary arcs, say
$a\bra{\alpha} a\bra{\rt\alpha}$.
The difference on the left hand side of (\ref{eq:hold})
is equal to
$$H_\beta(z)-H_\beta(z')=
\Prob\br{Q_\beta(z)\setminus Q_\beta(z')}-
\Prob\br{Q_\beta(z')\setminus Q_\beta(z)}~,$$
so it is enough to estimate the terms
on the right hand side of the above equation.
But for either of the corresponding events to occur,
the interval $[zz']$ should be joint to the arcs
$a\bra{\rt\beta} a\bra{\rt^2\beta}$,
$a\bra{\beta} a\bra{\rt^2\beta}$, and $a\bra{\beta} a\bra{\rt\beta}$,
by the yellow, blue, and blue clusters correspondingly
-- the reasoning is similar to that of the Lemma~\ref{lem:rot}.
Particularly, there is a monochrome cluster
connecting  $[zz']$ to $a\bra{\alpha} a\bra{\rt\alpha}$.

But $[zz']$ can be separated from $a\bra{\alpha} a\bra{\rt\alpha}$
by $|\log|z-z'||-c$ disjoint discrete annuli of fixed moduli.
By the Lemma~\ref{lem:rsw}
probability of the existence of a monochrome cluster traversing
such an annulus is bounded from above by some $q<1$,
regardless of size.
So the occurrence of a monochrome
cluster connecting $[zz']$ to $a\bra{\un} a\bra{\rt}$,
implies a simultaneous occurrence of
$|\log|z-z'||-c$ independent events of probability at most $q<1$,
and hence we infer
$$\abs{H_\beta(z)-H_\beta(z')}\le2q^{|\log|z-z'||-c}
=2q^{-c}\,\abs{z-z'}^{|\log q|}~,$$
and the desired estimate (\ref{eq:hold}) follows.

Note that as a particular consequence of this estimate
for the case $z'=z+\eta$ we obtain
\begin{equation}\label{eq:hol}
P_\beta(z,\eta)~\le~C'\mesh^{\varepsilon}~.
\end{equation}

To prove the Lemma it remains to check the boundary values of $H_\beta$.
Points on the arc $a\bra{\rt\beta}a\bra{\rt^2\beta}$
cannot be separated from it by a blue path,
so boundary values of $H_\beta$ on this arc are identically zero.
The boundary value of $H_\beta$ at the point $a\bra{\beta}$ 
tends to $1$ as $\mesh\to0$ also because of the Lemma~\ref{lem:rsw}.
Indeed, there are $|\log\mesh|-c$ disjoint discrete annuli of fixed shape
around $a\bra{\beta}$, each of them contains a blue cluster 
circumventing it independently with probability bounded from below by 
some $p>0$. Thus probability of $a\bra{\beta}$ being separated
from the arc $a\bra{\rt\beta}a\bra{\rt^2\beta}$ is
at least $1-(1-p)^{|\log\mesh|-c}$, which tends to $1$ as $\mesh\to0$.
This concludes the Proof.
\end{proof}

Take some equilateral triangular contour $\Gamma$
with vertices in the centers of the lattice triangles and
with bottom side parallel to the real axis.
Denote vertices of $\Gamma$ by $x\bra{\un},x\bra{\rt},x\bra{\rt^2}$ counterclockwise,
starting with the top one.
For a function $H(z)$ 
we define the discrete contour integral 
by
$$\doint_\Gamma H(z)dz~:=~
\mesh\sum_{z\in x\bra{\rt}x\bra{\rt^2}}H(z)+
\mesh\rt\sum_{z\in x\bra{\rt^2}x\bra{\un}}H(z)+
\mesh\rt^2\sum_{z\in x\bra{\un}x\bra{\rt}}H(z)~,$$
where sums are taken over centers $z$ of lattice triangles,
lying in the corresponding intervals.

\begin{lem}[Contour integrals vanish]\label{lem:md}
For any equilateral triangular contour $\Gamma\subset\Omega$ of length $\ell$
with vertices in the centers of the lattice triangles,
with bottom side parallel to the real axis,
and any $\beta\in\brs{1,\rt,\rt^2}$
one has
$$\doint_\Gamma H^\mesh_\beta(z)\,dz~=
~\doint_\Gamma \frac1{\rt}H^\mesh_{\rt\beta}(z)\,dz~+~
O(\ell\,\mesh^{\varepsilon})~.$$
\end{lem}

\begin{proof}
Color all triangles in chess-board fashion, so that
triangles with centers on $\Gamma$ are colored black.
Denote by ${\cal B}$ the set of all centers of black triangles,
lying on or inside $\Gamma$,
and by ${\cal W}$ the set of all centers of white triangles,
lying inside $\Gamma$.

Fix some $\alpha\in\brs{1,\rt,\rt^2}$.
Take $\eta$ to be of length $\mesh/\sqrt{3}$
collinear with $e^{\pi i/6}\,\br{x\bra{\rt^2\alpha}-x\bra{\rt\alpha}}$  
and denote $\eta':=e^{\pi i/3}\eta$.
We can write using (\ref{eq:der}) and (\ref{eq:rot}):
\begin{equation}\label{eq:ch1}
\begin{aligned}
\sum_{z\in{\cal B}\setminus x\bra{\alpha}x\bra{\rt^2\alpha}}& 
\br{H_\beta(z+\eta)-H_\beta(z)}
~\stackrel{(\ref{eq:der})}{=} 
~\sum_{z\in{\cal B}\setminus x\bra{\alpha}x\bra{\rt^2\alpha}} 
\br{P_\beta(z,\eta)-P_\beta(z+\eta,-\eta)}\\
&\stackrel{(\ref{eq:rot})}{=} 
\sum_{z\in{\cal B}\setminus x\bra{\alpha}x\bra{\rt^2\alpha}} 
\br{P_{\rt\beta}(z,\rt\eta)-P_{\rt\beta}(z+\eta,-\rt\eta)}\\
&\stackrel{(*)}{=}
\sum_{z\in{\cal B}\setminus x\bra{\alpha}x\bra{\rt\alpha}} 
\br{P_{\rt\beta}(z,\rt\eta)-P_{\rt\beta}(z+\rt\eta,-\rt\eta)}
~+~O(\ell\mesh^{\varepsilon-1})\\
&\stackrel{(\ref{eq:der})}{=}
\sum_{z\in{\cal B}\setminus x\bra{\alpha}x\bra{\rt\alpha}} 
\br{H_{\rt\beta}(z+\rt\eta)-H_{\rt\beta}(z)}
~+~O(\ell\mesh^{\varepsilon-1})~.
\end{aligned}
\end{equation}
In the identity (*) we also used that two sides differ by at most $O(\ell\mesh^{-1})$ 
terms (number of vertices on $\Gamma$),
which are of the order $O(\mesh^{\varepsilon})$ by(\ref{eq:hol}).
Similarly one shows that
\begin{equation}\label{eq:ch2}
\sum_{z\in{\cal W}}\br{H_\beta(z+\eta')-H_\beta(z)}
~=~\sum_{z\in{\cal W}} 
\br{H_{\rt\beta}(z+\rt\eta')-H_{\rt\beta}(z)}
~+~O(\ell\mesh^{\varepsilon-1})~.
\end{equation}

Now, combining (\ref{eq:ch1}) and  (\ref{eq:ch2})
we can write (using ``telescoping sums'')
\begin{eqnarray*}
&{\displaystyle\sum_{z\in x\bra{\alpha}x\bra{\rt^2\alpha}}}&H_\beta(z)-
\sum_{z\in x\bra{\rt\alpha}x\bra{\rt^2\alpha}} H_\beta(z)=\\
&=&\sum_{z\in{\cal B}\setminus x\bra{\alpha}x\bra{\rt^2\alpha}} 
\br{H_\beta(z+\eta)-H_\beta(z)}+
\sum_{z\in{\cal W}}\br{H_\beta(z+\rt\eta')-H_\beta(z)}\\
&=&\sum_{z\in{\cal B}\setminus x\bra{\alpha}x\bra{\rt\alpha}} 
\br{H_{\rt\beta}(z+\rt\eta)-H_{\rt\beta}(z)}
+\sum_{z\in{\cal W}} 
\br{H_{\rt\beta}(z+\rt\eta')-H_{\rt\beta}(z)}
+O(\ell\mesh^{\varepsilon-1})\\
&=&\sum_{z\in x\bra{\alpha}x\bra{\rt\alpha}}H_{\rt\beta}(z)-
\sum_{z\in x\bra{\alpha}x\bra{\rt^2\alpha}}H_{\rt\beta}(z)
~+~O(\ell\mesh^{\varepsilon-1})~.
\end{eqnarray*}
So for any $\alpha\in\brs{1,\rt,\rt^2}$
one has
\begin{equation}\label{eq:dm}
\sum_{z\in x\bra{\alpha}x\bra{\rt^2\alpha}} H_\beta(z)-
\sum_{z\in x\bra{\rt\alpha}x\bra{\rt^2\alpha}} H_\beta(z)
=\sum_{z\in x\bra{\alpha}x\bra{\rt\alpha}} H_{\rt\beta}(z)-
\sum_{z\in x\bra{\alpha}x\bra{\rt^2\alpha}} H_{\rt\beta}(z)
+O(\ell\mesh^{\varepsilon-1})~.
\end{equation}
Summing three copies of equation~(\ref{eq:dm}) above with coefficients and
different values of $\alpha$ plugged in:
$$-\frac{\mesh}{2}{\mathrm(equation~\ref{eq:dm})}|_{\alpha=1}
-i\frac{\mesh\sqrt3}{2}{\mathrm(equation~\ref{eq:dm})}|_{\alpha=\rt}
+\frac{\mesh}{2}{\mathrm(equation~\ref{eq:dm})}|_{\alpha=\rt^2}~,$$
we arrive at the desired identity and prove the Lemma.
\end{proof}

By Lemma~\ref{lem:hold}
the H\"older norms of functions $H^\mesh_\alpha$, $\alpha\in\brs{1,\rt,\rt^2}$
are uniformly bounded,
hence from any sequence of such functions with $\mesh\to0$
one can chose a uniformly converging subsequence.
Therefore to show that the functions $H^\mesh_\alpha$ 
converge uniformly, as $\mesh\to0$,
to the functions $h_\alpha$ and prove Theorem~\ref{thm:triangle},
it is sufficient prove the following
\begin{lem}\label{lem:lim}
Assume that for some subsequence $\mesh_j\to0$
the functions $H^{\mesh_j}_\alpha$ converge uniformly in $\Omega$
to some functions $f_\alpha$.
Then $f_\alpha\equiv h_\alpha$.
\end{lem}

\begin{proof}
Clearly the discrete contour integrals $\doint H^{\mesh_j}_\beta$
converge to the usual contour integrals $\oint f_\beta$.
Then using Lemma~\ref{lem:md} 
we conclude
(one may need to shift contour by $\le\mesh$ so that
its vertices are in the centers of lattice triangles,
but this perturbation does not affect the limit), 
that for any equilateral triangular 
contour $\Gamma\subset\Omega$ of length $\ell$,
with bottom side parallel to the real axis,
and any $\beta\in\brs{1,\rt,\rt^2}$
one has
\begin{equation}\label{eq:mor1}
\oint_\Gamma f_\beta(z)\,dz~=
~\oint_\Gamma \frac1{\rt}f_{\rt\beta}(z)\,dz~.
\end{equation}
Take any $\alpha\in\brs{1,\rt,\rt^2}$.
Adding two copies of (\ref{eq:mor1}) with coefficients and
different values of $\beta$ plugged in:
$${\mathrm (equation~\ref{eq:mor1})}|_{\beta=\alpha}-\br{\frac12+\frac{i}{2\sqrt3}}
\cdot{\mathrm (equation~\ref{eq:mor1})}|_{\beta=\rt\alpha}~,$$
we obtain
\begin{equation}\label{eq:mor2}
\oint_\Gamma \br{f_\alpha(z)\,+\,\frac{i}{\sqrt3}
\br{f_{\rt\alpha}(z)-f_{\rt^2\alpha}(z)}}\,dz~=~0~.
\end{equation}
By Morera's theorem we deduce that
the function in (\ref{eq:mor2}) is analytic,
and hence for any $\alpha\in\brs{1,\rt,\rt^2}$
the function $f_\alpha$ is harmonic
with harmonic conjugate
$\frac{1}{\sqrt3}\br{f_{\rt\alpha}(z)-f_{\rt^2\alpha}(z)}$.
It follows that for any $\alpha\in\brs{1,\rt,\rt^2}$
and any unit vector $\eta$
\begin{equation}\label{eq:cr}
\frac{\partial}{\partial\eta}f_\alpha~=~
\frac{\partial}{\partial(\rt\eta)}f_{\rt\alpha}~.
\end{equation}
Hence
$$\begin{aligned}
\frac{\partial}{\partial(\rt\tang)}f_{\alpha}\equiv
\frac{\partial}{\partial\tang}f_{\rt^2\alpha}\equiv
\frac{\partial}{\partial\tang}0\equiv0&
{\mathrm~~for~points~on~the~arc~}a\bra{\alpha}a\bra{\rt\alpha}~,\\
\frac{\partial}{\partial(-\rt^2\tang)}f_{\alpha}\equiv
\frac{\partial}{\partial(-\tang)}f_{\rt\alpha}\equiv
\frac{\partial}{\partial(-\tang)}0\equiv0&
{\mathrm~~for~points~on~the~arc~}a\bra{\rt^2\alpha}a\bra{\alpha}~.
\end{aligned}$$
The identities above should be understood as
holding in the limit as $z\in\Omega$ tends to the corresponding boundary arc.
Also by the Lemma~\ref{lem:hold}
the boundary values of
$f_\alpha$ are equal to zero on the arc $a\br{\rt\alpha}a\br{\rt^2\alpha}$
and to $1$ at the point $a\br{\alpha}$.

Summing it up, we conclude
that $f_\alpha$'s satisfy the mixed Dirichlet-Neumann problem (\ref{eq:ndp}),
which has a unique solution.
Thus $f_\alpha=h_\alpha$, and we have proven
the Lemma and the Theorem.
\end{proof}

\section{Continuum scaling limits}\label{sec:scaling}

Consider some domain $\Omega$ with three points (or prime ends) 
$a$, $b$, and $c$ on the boundary, named counterclockwise.
For a discrete percolation configuration we can define 
``the lowest blue crossing,''
$\curlow=\curlow_\mesh=\curlow_\bl(a,b,c)$, as
the closest to the arc $ab$,
simple blue curve going from  the arc $ca$ to the arc $bc$.
For simplicity, if no such crossing exists,
we assume it to consist of one vertex, closest to $c$.
This does not influence our reasoning since by the Lemma~\ref{lem:rsw}
the probability of the existence of a genuine crossing
tends to $1$ as $\mesh\to0$.
Similarly, we define ``the highest yellow crossing,''
$\curlow_\ye(c,a,b)$, as the
closest to the arc $ca$,
simple yellow curve going from  the arc $bc$ to the arc $ab$.

These curves end at neighboring vertices on the arc $bc$,
so we can join them to obtain a simple curve
$\curhull=\curhull_\mesh=\curhull_{\ye,\bl}(a,b,c)$, which starts on the arc
$ab$ and ends on the arc $ca$, touching the arc $bc$ in the middle.
The law of $\curhull_\mesh$ is a probability measure $\lawhull_\mesh$,
supported on the space $\holcur$ of H\"older curves
(we identify a discrete curve on the $\mesh$-lattice with the
corresponding broken line).

\begin{thm}[Outer boundary]\label{thm:hull}
As $\mesh\to0$, the law of $\curhull$
converges to a law $\lawhull$
on H\"older self-avoiding paths from $ab$ to $ac$,
touching $bc$ at one point, which separates
them into two simple ``halves.''
The law $\lawhull$ is a conformal invariant of the configuration $(\Omega,a,b,c)$.
\end{thm}

Note that in the scaling limit two halves of the $\curhull$
are the outer boundaries (inside $\Omega$) 
of all yellow percolation clusters intersecting the arc $ab$ 
and viewed from $c$,
and all blue percolation clusters intersecting the arc $ca$
and viewed from $b$.



In the proof we automatically obtain, that
the law of $\curhull$ is CCI (Completely Conformally Invariant
in the language of \qcite{Lawler-Werner-universality,Werner-ecm}),
and points where $\curhull$ touches the three arcs satisfy the Cardy-Carleson law
(Remark~\ref{rem:hull} also gives that $h_{a}(c)$ is the probability of a point
to be separated by $\curhull$ from $b$ and $c$).
These two properties determine the law uniquely (cf. \qcite{Werner-ecm}).
Hence using \qcite{Lawler-Schramm-Werner-i,Werner-ecm} we conclude that
the law of $\curhull$ coincides with two other laws, enjoying the same properties:
with the laws of the hull of Oded Schramm's chordal SLE${}_6$
(started at $a$, aimed at some point in $bc$ and stopped upon hitting $bc$)
or hull of the reflected Brownian motion (started at $a$, reflected on $ab$ and $ac$
at $\frac\pi3$-angle pointing towards $bc$,
and stopped upon hitting $bc$), 
which was considered by Wendelin Werner in \qcite{Werner-ecm}.

\begin{cor}\label{cor:hull}
The law $\lawhull$ coincides with the law of the boundary of the hull
of chordal SLE${}_6$ (or reflected Brownian motion) started at $a$
and stopped upon hitting $bc$.
\end{cor}

By the work 
\qcite{Lawler-Schramm-Werner-43}
of Greg Lawler, Oded Schramm, and Wendelin Werner
we know the dimension of the outer boundary of Brownian motion, so we deduce 

\begin{cor}
The curves $\curhull$ and $\curlow$ 
(and hence the outer boundary of a percolation cluster)
have Hausdorff dimension $\frac43$ almost surely.
\end{cor}

\begin{rem}
Lowest crossing depends only on the area ``below'' it,
and all our considerations are stable under perturbations of domains,
since Cardy's formula is.
So one can continue working with the area ``above'' the lowest
crossing, finding crossings for some configurations there.
Using an inductive procedure it is not difficult to construct  in this way
scaling limits for laminations and backbones
(and then in a similar inductive way pass to the full percolation  configuration).
\end{rem}

Consider now a domain $\Omega$ with two boundary points $a$ and $b$.
For any percolation  configuration there is a unique curve $\curexp$
along the edges of the dual hexagonal lattice,
which goes from $a$ to $b$ separating
the blue clusters intersecting the arc $ab$
from the yellow clusters intersecting the arc $ba$.
This curve is called the ``exploration process''
(one can actually introduce time)
or ``external perimeter''
(of all blue clusters intersecting $ab$, as viewed from $ba$).

We say that a path is {\it self-avoiding}
if it does not have ``transversal self intersections''
(maybe a more appropriate term would
be ``non-self-traversing'').

\begin{thm}[External perimeter]\label{thm:explor}
As $\mesh\to0$, the law of the exploration process $\curexp$
converges to a law $\lawexp$
on H\"older self-avoiding paths from $a$ to $b$.
The law is a conformal invariant of the configuration $(\Omega,a,b)$.
\end{thm}

It is possible to prove (using known crossing exponent for ``$5$ arms'')
that the curve $\curexp$ is almost surely not simple.
In proving that the subsequential limit of the laws $\lawexp_\mesh$
is uniquely determined, we only use its properties 
(locality and Cardy's formula)
valid also for Oded Schramm's SLE${}_6$ process
by \qcite{Lawler-Schramm-Werner-i,Schramm-lerw},
and so we arrive at
\begin{cor}\label{cor:sle}
The law $\lawexp$ coincides
with that of the Schramm's chordal SLE${}_6$
started at $a$ and aiming at $b$.
\end{cor}

As discussed above, the discrete percolation  configuration
in the whole plane can be represented
by a collection of all external perimeters:
(as a collection of nested pairwise disjoint oriented simple loops). 
It is not difficult (though a bit technical) to obtain from the results above 
the following
\begin{thm}[Collection of external perimeters]\label{thm:colper}
As $\mesh\to0$, the law of the collection
of all external perimeters converges to a
law $\lawcolper$ on collections of mutually- and self-avoiding
oriented nested H\"older loops.
This law (and its restriction to any domain)
is conformally invariant.
\end{thm}
We sketch the proof below. In a subsequent paper we intend
to give a different, perhaps more conceptual, proof of this Theorem.

\subsection{Proof for the outer boundary} 

By the Lemma~\ref{lem:holnorm} (which comes from 
the Theorem A.1 in 
\qcite{Aizenman-Burchard} by M.~Aizenman and A.~Burchard)
the family of measures  $\brs{\lawhull_\mesh}_\mesh$ 
is weakly precompact in $\holcur$.
Indeed, $\lim_{M\to\infty}\lawhull_\mesh(\holcur\setminus\holcurb)=0$
uniformly in $\mesh$, whereas the sets $\holcurb$ are compact.

Thus one can chose a sequence $\mesh_j\to0$ so that
the subsequence $\brs{\lawhull_{\mesh_j}}$
converges weakly to a measure $\lawhull$.
It is sufficient to show that the latter measure
is independent of the chosen subsequence.
For the rest of the section we fix the sequence $\{\mesh_j\}$.

Since the property of being self-avoiding is preserved in the limit,
the measure $\lawhull$ is supported on self-avoiding H\"older curves,
starting on the arc $ab$, touching the arc $bc$
and ending on the arc $ac$.
The two ``halves'' of the curve $\curhull$ will almost surely be simple
(However, the whole curve need not be).
Indeed, the probability of a simple curve becoming non-simple 
in the scaling limit is zero:
it would imply six ``arms,'' which 
happens with probability zero by Lemma~\ref{lem:plane6}.
Similarly one shows, that almost surely
the curve $\curhull$ touches the arc $bc$ at a unique point
(otherwise there would be three ``arms'' going
to a boundary point, which happens with probability zero).
Denote by $\holcurhull$ the set of the H\"older curves
satisfying the properties above.

For two disjoint simple curves 
$\eta$ (going from $p\in ab$ to $q\in bc$ and denote $pq$ below)
and $\eta'$ (going from $r\in bc$ to $s\in ac$ and denoted $rs$ below),
let $A_{\eta,\eta'}$ be the event of $\curhull$ 
lying completely inside the conformal pentagon $apqrs$.
It is easy to see that the event $A_{\eta,\eta'}$ is equivalent
to the a yellow crossing from $ap$ to $pqrs$ within $apqrs$,
which is closest to $as$, ending inside $qr$.
Since this probability is given by the Cardy-Carleson law,
it is uniquely determined and so is the probability of $A_{\eta,\eta'}$.

It remains to notice that the events $A_{\eta,\eta'}$ generate 
(by disjoint unions and complements)
the restriction of the $\sigma$-algebra $\siglaws$ to
the set $\holcurhull$, which supports the measure $\lawhull$.
Hence the measure $\lawhull$ is independent of the subsequence chosen.

Note that we automatically obtain that $\lawhull$ is conformally
invariant, since so is the Cardy-Carleson law.
The Corollaries follow, since the properties we used are
satisfied by the hull of reflected Brownian motion (see \qcite{Werner-ecm})
and by the hull of the Schramm's SLE$_{6}$ (see \qcite{Lawler-Schramm-Werner-i}).

\subsection{Proof for the external perimeter} 

By the Lemma~\ref{lem:holnorm} 
the family of measures  $\brs{\lawexp_\mesh}_\mesh$ 
is weakly precompact in $\holcur$.
Thus one can choose a sequence $\mesh_j\to0$ so that
the subsequence $\brs{\lawexp_{\mesh_j}}$
converges weakly to a measure $\lawexp$.
Again, it is sufficient to show that the latter measure
is independent of the chosen subsequence.
For the rest of the section we fix the sequence $\{\mesh_j\}$.

Since such properties are preserved in the limit,
the measure $\lawexp$ is supported on self-avoiding H\"older curves,
starting at $a$, and terminating at $b$.
When $\curexp(t)$ is some parametrization of such a curve by the interval $[0,1]$,
denote by $\hull(t)$ its hull at time $t$,
i.e. the closure of $\Omega\setminus\Omega(t)$.
Here $\Omega(t)$ denotes the component of connectivity
of $\Omega\setminus\curexp[0,t]$, containing the point $b$.
Then the hull grows, i.e. it is strictly larger for larger values of $t$, 
meaning that the endpoint of $\curexp[0,t]$ is always ``visible'' from $b$.
Indeed, otherwise the external perimeter would enter
a zero width fjord, implying eight ``arms,'' which happens 
with probability zero by Lemma~\ref{lem:plane6}.

Denote by $\holcurexp$ the set of the H\"older curves
satisfying the properties above.
It is easier to work with perimeters by parameterizing them (e.g. by diameter).
One can do it in the following way.
Namely, to each $\gamma\in\holcurexp$ and $\epsilon>0$
we associate a broken line $\gamma_\epsilon=\phi_\epsilon(\gamma)$ 
by the following inductive procedure.
The first point is $\gamma^0_\epsilon:=\gamma(t_0)=a$.
Then $\gamma^{j+1}$ is the first exit of $\gamma[t_j,1]$
from the ball $B(\gamma^{j},\epsilon)$.
If it never exits this ball, we set $\gamma^{j+1}:=b$
and end the procedure.

The law of curves $\gamma_\epsilon$ is uniquely determined.
In fact it is easy to see by induction that the laws
of $t_j$, $\hull(t_j)$, and hence $\gamma^{j}$ are uniquely determined: 
in fact, the law of $\hull(t_{j+1})\setminus\hull(t_j)$
coincides with the law $\lawhull$ inside $B(\gamma_j,\epsilon)\setminus\hull(t_j)$.
We use that percolation is local,
and the hull depends only on percolation inside it
(stability of Cardy's formula under perturbations
of the boundary is also used here).

Thus the measure $\phi_{\epsilon}^{-1}(\mu)$ is uniquely determined.
Recall that $\lim_{M\to\infty}\lawexp(\holcur\setminus\holcurb)=0$
and also note that $\gamma_\epsilon$ converges to $\gamma$ uniformly in $\holcurb$.
Therefore, $\phi_{\epsilon}^{-1}(\mu)$ converge weakly to $\mu$,
and $\mu$ is uniquely determined.

As before, we automatically obtain that $\lawexp$ is conformally
invariant, since so are the Cardy-Carleson law and $\lawhull$.
The Corollary follows, since the properties we used are
satisfied by the Schramm's SLE$_{6}$ (see \qcite{Lawler-Schramm-Werner-i}).



\subsection{Sketch of the proof for the full limit}\label{subsec:cc}

Again, \qcite{Aizenman-Burchard} gives enough information
to obtain a subsequential scaling limit: a measure $\lawcolper$ on the mentioned space,
and it is sufficient to prove that this measure is uniquely determined.
One of a few possible ways is to proceed as follows.

For any domain $\Omega$ with two boundary points 
$a$ and $b$
to the law $\lawcolper$ corresponds the law $\lawexp=\lawexp_\Omega$
of external perimeter in $\Omega$, which is independent from the choice of 
a subsequence determining $\lawcolper$. 
For any particular external perimeter curve
$\curexp\in\holcurexp(\Omega)$, consider
the connected components $\{\Omega'\}$ of $\Omega\setminus\curexp$.
Usual ``number of arms'' considerations
show that almost surely no other external perimeters can go from
one such component to another within $\Omega$
(one has to also use even more classical estimate
for $3$ arms in a half-plane).
Also the law $\lawexp$ depends only on $\curexp$
(or rather on its ``infinitesimal neighborhood'').
So for any component $\Omega'$ bordering the boundary on the arc $a'b'$
we can do the similar considerations, retrieving uniquely the law
of the external perimeter within $\Omega'$
from $a'$ to $b'$ (or in the opposite direction,
depending on the location of $a'$ and $b'$: 
whether it is $ba$ or $ab$).
There is a subtle point here: one has to check that in some
(related to conformal geometry) sense law of
external perimeter is stable under perturbations of a domain.
Continuing by induction,
 one retrieves (in the limit) uniquely the law
of all arcs of curves (from the full collection of external perimeters)
with endpoints on the boundary of $\Omega$.
Now performing a new inductive procedure with this law, one
can obtain the full scaling limit.

\section{Appendix: Technical estimates}\label{sec:tech}

We gathered in this Section some known
(see \qcite{Aizenman-statphys,Aizenman-Burchard,Grimmett-book99}) estimates we need.

\begin{lem}[Russo-Seymour-Welsh estimates]\label{lem:rsw}
For a conformal rectangle (or an annulus) of a fixed shape
probability of the existence of a monochrome crossing
is contained in an interval $[p,q]$.
The constants $p>0$, $q<1$ depend on the shape, but 
not on the size of the rectangle or mesh of the lattice.
\end{lem}
This Lemma is needed to establish H\"older continuity of harmonic conformal invariants.
For bond percolation such a statement can be found 
in 11.70 of G.~Grimmett's book \qcite{Grimmett-book99}.
The proof is based on self-duality property 
(forcing the crossing probabilities of symmetric shapes to be $\frac12$),
so it equally applies to critical site percolation on triangular lattice.

\begin{lem}\label{lem:holnorm}
For discrete percolation all the non-repeating paths supported on 
the connected clusters (or external perimeters of clusters)
within some compact region
can be simultaneously parameterized by H\"older
continuous functions on $[0,1]$,
whose H\"older norms are (uniformly in the mesh $\mesh$)
stochastically bounded.
\end{lem}
This Lemma is used to establish precompactness
of laws describing discrete percolations.
It is stated as Theorem A.1 at the end of \qcite{Aizenman-Burchard}
by M.~Aizenman and A.~Burchard 
for non-repeating paths, and the main theorems of the paper clearly
apply to external perimeters as well.



\begin{lem}\label{lem:plane6}
For two concentric balls of radii $r$ and $R$
probability of existence for discrete percolation
of $5$ crossings (not all of the same color) between their boundaries
is $\le{\mathrm{const}}(r/R)^2$ uniformly in the mesh $\mesh$ of the lattice. 
Hence (e.g. by the Russo-Seymour-Welsh theory) the similar probability
for $6$ crossings is $\le{\mathrm{const}}(r/R)^{2+\epsilon}$,
and it immediately follows that in the scaling limit
almost surely there is no point with $6$ curves
(not all of the same color) incoming.
\end{lem}
We use this Lemma (see \qcite{Aizenman-statphys} by M.~Aizenman)
to show that simple curves cannot ``collapse'' when we
pass to the subsequential limit.
For general number of arms, the exact exponents were predicted
by M.~Aizenman, B.~Duplantier, and H.~Aharony in
\qcite{Aizenman-Duplantier-Aharony}
(see also related papers of B.~Duplantier,
e.g. \qcite{Duplantier-harmonic} and references therein).
Up to now there was no direct proof, but
since exponents are known for SLE${}_6$,
it is not difficult to write a proof on the basis
of the current paper and work of G.~Lawler, O.~Schramm, and W.~Werner
\qcite{Lawler-Schramm-Werner-i,Lawler-Schramm-Werner-ii,Lawler-Schramm-Werner-iii}.

We only need a very special case of $5$ arms 
which can be calculated directly
(along with $2$ or $3$ arms in the half-plane,
since they have exponents $1$ or $2$ and some
``derivative representation.'')

\end{document}